\documentclass[final,leqno]{siamltex}
\usepackage{moreverb}
\pdfoutput=1
\usepackage{graphicx,subfigure,wrapfig}
\usepackage{algorithm}
\usepackage{algpseudocode}
\usepackage{color}
\usepackage{amsmath}
\usepackage{txfonts}
\usepackage{afterpage}
\usepackage[margin=1in]{geometry}
\usepackage{mathrsfs}
\usepackage[american]{circuitikz}
\usepackage{mathrsfs}
\usepackage{pstricks}
\usepackage{makecell}
\usepackage{threeparttable}
\usepackage{booktabs}

\newtheorem{mdefinition}{Definition}
\newtheorem{massumption}{Assumption}

\newtheorem{mtheorem}{Theorem}
\newcommand{\ls}{\color{black}}
\renewcommand{\mathbf}{\boldsymbol}

\newcommand{\mb}{\mathbf}

\title{A hierarchical neural hybrid method for failure probability estimation}

\author{Ke Li$^{\displaystyle 1}$ \thanks{School of Information Science and Technology, ShanghaiTech University, Shanghai, 200120, China. $^{\displaystyle 1}$ Equal contributions. 
		({\tt \{like1, tangkj, liaoqf\}@shanghaitech.edu.cn}).}
        \and Kejun Tang$^{\displaystyle 1}$ \footnotemark[1] 
         \and Jinglai Li\thanks{Department of Mathematical Sciences, University of Liverpool, Liverpool, UK,
        ({\tt Jinglai.Li@liverpool.ac.uk}).}
        \and Tianfan Wu\thanks{Viterbi School of Engineering, University of Southern California, Los Angeles, USA,
        ({\tt tianfanw@usc.edu}).}  
    	\and Qifeng Liao\footnotemark[1] \thanks{Corresponding author.}	
    	}
    
\begin{document}

\maketitle

\begin{abstract}
Failure probability evaluation for complex physical and engineering systems governed by partial differential equations (PDEs)
are computationally intensive, especially when high-dimensional random parameters are involved. 
Since standard numerical schemes for solving these complex PDEs are expensive,
traditional Monte Carlo methods which require repeatedly solving PDEs are infeasible. 
Alternative approaches which are typically the surrogate based methods suffer from the so-called ``curse of dimensionality'',
which limits their application to problems with high-dimensional parameters.
For this purpose, we develop a novel hierarchical neural hybrid (HNH) method to efficiently compute failure probabilities of these challenging high-dimensional  problems.
Especially,  multifidelity surrogates are constructed based on neural networks with different levels of layers, 
such that expensive highfidelity surrogates are adapted only when the parameters are in the suspicious domain.
The efficiency of our new HNH method is theoretically analyzed and is demonstrated with numerical experiments.
From numerical results, we show that to achieve an accuracy in estimating the rare failure probability (e.g., $10^{-5}$), the traditional Monte Carlo 
method needs to solve PDEs more than a  million times, while our HNH only requires solving them a few thousand times.
\end{abstract}


\pagestyle{myheadings}
\thispagestyle{plain}
\markboth{K. LI, K. TANG, J. Li, T. Wu, AND Q. Liao}{Hierarchical Neural Hybrid Method}

\section{Introduction}
% Section 1: introduction and related work
Due to lack of knowledge or measurement of realistic model parameters, modern complex physical and engineering systems are often modeled by partial differential equations (PDEs) 
with high-dimensional random parameters. For example, groundwater flow problems are modeled by stochastic diffusion equations, and acoustic scattering problems are modeled 
by Helmholtz equations with random inputs. 
When conducting risk management, it is essential to compute failure probabilities of these stochastic PDE models.  
A standard method to compute the failure probability is the traditional  Monte Carlo sampling method \cite{metropolis1949monte}.
However, this method requires repeatedly solving complex PDEs to generate a large number samples to capture failure probabilities associated with rare events.
There are two main computational challenges: first, it is expensive to solve each complex PDE using standard numerical schemes (e.g., finite elements \cite{elman2014finite});
second, the complex PDEs need to be repeatedly solved many times for computing  failure probabilities.

%Gradient-based method and simulation-based method are the mainly two types of approaches to estimate failure probability. First-order reliability method (FORM) \cite{der1998multiple,zhao1999general}, second-order reliability method (SORM) \cite{hohenbichler1987new,der1987second,kiureghian1991efficient} and a statistical technical response surface method (RSM) \cite{rajashekhar1993new,bucher1990fast,kaymaz2005response} are classical gradient-based methods. Simulation-based method relies on Monte Carlo sampling (MCS) and the failure probability can be approximated by the ratio between number of failure samples and the total amount. Using standard MC method means generating a large quantity of random variables from the given distribution and solving the complex system for each set of variable. Conditional simulation \cite{ayyub1995simulation}, importance sampling (IS) \cite{melchers1989importance} and sequential Monte Carlo (SMC) \cite{Bouchard2012Phylogenetic,Bouchard2014Sequential} are efficient sampling improvements. 

As alternatives, gradient-based methods and simulation-based methods are developed and are briefly reviewed as follows. 
Classical gradient-based methods are the first-order reliability method (FORM) \cite{der1998multiple}, 
second-order reliability method (SORM) \cite{hohenbichler1987new} and a statistical technical response surface method (RSM) \cite{rajashekhar1993new,bucher1990fast}. 
Moreover,  Simulation-based method relies on Monte Carlo (MC) sampling and the failure probability can be approximated by the ratio between number of failure samples and the total amount. 
To deal with rare events, conditional simulation \cite{ayyub1995simulation} and importance sampling (IS) \cite{melchers1989importance} are efficient improvements of the traditional Monte Carlo method.

%In this paper, we propose a neural hybrid method and a hierarchical neural hybrid method to solve former challenges. Our methods combine the advantage of MC, RSM and neural network to solve high-dimensional failure problems, which are expensive or infeasible to solve by normal approaches. Hierarchical neural hybrid (HNH) method is designed with following properties. First, considering that repeatedly solving a fine-fidelity surrogate is still costly, HNH method constructs multifidelity surrogates to further speed up without sacrificing accuracy. Second, only a few samples are in the suspicious domain, it requires the fine-fidelity model and the total cost will be much smaller than before. Third, our method is a universal solver with neural networks since it can approximate any complex system with an acceptable accuracy \cite{csaji2001approximation}. The efficiency of our new HNH method is theoretically analyzed and it is mathematically proofed that our method is convergent and unbiased. The results of numerical experiments demonstrate that our methods are highly efficient algorithm with remarkable accuracy.%the We construct a surrogate by using 

However, the above methods suffer from the ``curse of dimensionality'', which limits their applications for complex problems with high-dimensional inputs.  
	
%\hl{perform well via reduced basis ANOVA decomposition} \cite{liao2019adaptive}\hl{ and Gaussian process} \cite{rasmussen2003gaussian}.
In this work, we propose a novel state-of-the-art hierarchical neural hybrid (HNH) method to resolve the challenges discussed above.
Our method combines the advantages of MC, RSM and neural network to solve high-dimensional failure probability problems.
The main advantages of HNH are three-fold. 
First, HNH is a universal solver based on neural networks which can efficiently approximate any complex system including PDEs
with any given accuracy \cite{hornik1989multilayer}.
Second, multifidelity surrogates are constructed and expensive fine-fidelity surrogates are adopted only for samples close to the suspicious domain,
which results in an overall efficient computational procedure. 
%so that the  cost of the overall procedure is reduced without sacrificing accuracy.
Finally, only a few samples generated through solving given PDEs with standard numerical schemes are required to construct the surrogates and 
to modify the failure probability estimation.

%multilayer neural network to approximate the given complex system. We choose neural networks as the surrogate since it can approximate any complex system with an acceptable accuracy \cite{csaji2001approximation}. Then we employ hybrid method to combine sampling from surrogate and origin system. With hybrid method, only a few samples are in the suspicious domain, it requires the fine-fidelity model and the total cost will be much smaller than before. Moreover, these small part of samples decrease the loss of accuracy. Neural hybrid method is more efficient than direct MCS and offers a high accuracy. Considering that using a fine-fidelity surrogate is also costly, we propose a hierarchical neural hybrid (HNH) method which speeds up by using multifidelity models. The numerical results show that it is a highly efficient algorithm with remarkable accuracy.

\textbf{Our contributions are as follows.} 
\begin{itemize}
\item 
Our method is based on a combination of  neural networks and the hybrid method.
The hybrid method is an extremely effective approach which can capture failure probability using a few samples without losing accuracy. And our new method utilizes the feature of neural network as a universal approximation to solve complex systems with high-dimensional inputs. % As a universal approximator, we employ neural network to figure out high-dimensional problems which is infeasible to calculate by using traditional polynomial chaos. 

\item We propose a novel hierarchical neural hybrid (HNH) method. Given the fact  that a sufficient deep neural network surrogate is still expensive
to approximate complex systems, we employ multifidelity models to replace the single fine-fidelity deep model. Our method uses coarse-fidelity surrogates as a
 preconditioning scheme and uses fine-fidelity surrogates  to correct the estimation. 
 Our HNH method can significantly accelerate the computational procedure for failure probability estimation without sacrificing accuracy.
\end{itemize} 

\textbf{Paper organization.} The failure probability, the hybrid method and multilayer neural network are briefly reviewed in Section 2. 
Our novel hierarchical neural hybrid method is presented in Section 3,
where the rigorous error analysis for HNH is conducted. 
In Section 4, the efficiency of HNH is demonstrated with a high-dimensional structural safety problem, stochastic diffusion equations and Helmholtz equations. 
Finally some concluding remarks are offered in Section 5.

%\subsection{Related work}
%\textbf{Generalized polynomial chaos (gPC).} As an extension on polynomial chaos \cite{ghanem1991stochastic}, gPC is a widely used method for stochastic computation. A practical gPC called Wiener-Askey polynomial chaos from \cite{xiu2002wiener} can be rapidly constructed for high orders. In \cite{li2011efficient}, they construct a very high order gPC to approximate the given system. But in high-dimensional problems, the construction of gPC will be too difficult and expensive to guarantee the accuracy.

% \textbf{Rare event.} The probability of occurrence of given events is a extremely low value and the problem can be called as rare event. The RESTART method from \cite{villen1991restart} constructs a event $C$ and uses the combination of $P(C)$ and $P(A/C)$. It can be regarded as an application of $splitting$ \cite{hammersley1965monte}. Reverse-time models from \cite{frater1990fast,frater1988estimation} show a significantly greater speedup in using uncontrolled ALOHA system \cite{cottrell1983large} as an example. 

\section{Preliminaries }
\subsection{Failure probability}
In practical engineering problems, the various uncertainties ultimately affects the structural safety. Reliability analysis has become increasingly attracted in engineering analysis, and it can measure the structural safety by considering these uncertainties \cite{song2019new,huang2018stochastic}. In a general setting, let $Z$ be a $n_z-$dimensional random vector $Z=(Z_1,Z_2,\dots,Z_{n_z}):\Omega\to R^{n_z}$, where $F_Z(z) = \mathrm{Prob}(Z\leq z)$ is the distribution function, $\Omega$ is the probability space. Given a (scalar) limit state function $g(Z)$, $g(Z)<0$ defines a failure domain $\Omega_f$ and $g(Z)\geq 0$ defines safe domain. For some specific setting in section 4, we give additional definitions of $\Omega_{f}$. The failure probability $P_f$ is defined
\begin{equation} \label{fail_prob}
\small
P_f = \mathrm{Prob}(Z\in \Omega_f) = \int_{\Omega_f}^{}dF_Z(z) = \int\chi_{\Omega_f}(z)dF_Z(z),
\end{equation}
where $\chi$ means the characteristic function
\begin{equation}
\chi_{{\Omega_f}}(z) = \left\{
\begin{array}{rcl}
1 \ \ \mathrm{if} \ \ z\in { \Omega_f}, \\
0 \ \ \mathrm{if} \ \ z \notin {\Omega_f}.
\end{array}\right.
\end{equation} 

% Due to the complexity of function in real-world, it is usually infeasible to get the exact limit state function $g(Z)$, and inputs $Z$ are often random variables. So that we can regard $g(Z)$ as a relationship between random inputs and failure probability. Exact $g(z)$ can only be obtained by solving a complex stochastic system, so our idea is to construct a close surrogate with neural network which can overcome the weakness of gPC surrogate in high-dimension problems. For reducing the runtime using a high-fidelity model, we propose a method called hierarchical neural hybrid (HNH) method . HNH employs a hierarchy of surrogate models and modifies results from low-fidelity model by using high-fidelity model. 

%For the first benchmark in this paper, our goal is to minimize a deterministic cost function while the constraint is that failure probability is lower than a given threshold $\theta$. The mathematical representation can be frequently written as
%\begin{equation} \label{Cons}
%\begin{array}{lcl}
% \min \limits_{x\in D} f(x),\\
%\mathrm{s.t.} \ \ \ c(x) = \ln P(x) - \ln \theta \leq 0,
%\end{array}
%\end{equation}
%where x is the parameter and $\mathrm{D}$ is the design domain. In the next section, Monte Carlo simulation will be introduced as a basic but costly method to compute the probabilistic constrain.

\subsection{Hybrid method }
Hybrid method \cite{li2010evaluation} can enhance the performance of Monte Carlo sampling with surrogate models. Instead of computing (\ref{fail_prob}) directly, we employ surrogate model $\hat{g}$ to evaluate the failure probability
\begin{equation}\label{fail_surr}
\hat{P_f} = \int\chi_{\{\hat{g}(z)<0\}}(z)q(z)dz,
\end{equation}
where $q(z)$ denotes an arbitrary distribution of variable $z$. Specifically, we utilize MC to evaluate (\ref{fail_surr}) by generating samples $\{z^{(i)}\}_{i=1}^M$ from $q(z)$
\begin{equation}\label{fail_surr_MC}
\hat{P}^{mc}_f = \frac{1}{M}\sum_{i=1}^{M}\chi_{\{\hat{g}(z)<0\}}(z^{(i)}).
\end{equation}
where $M$ denotes the number of samples.

%With the idea of response surface method (RSM) \cite{au2001estimation,khuri2010response} and  design from \cite{li2011efficient}, we present a probability formulation associated with the state function $g$, surrogate model $\hat{g}$, and a real parameter $\gamma$ related to the accuracy of surrogate model. To compute the probability of $g(Z)<0$, we first evaluate the value of $\hat{g}(z)$. For a When $\hat{g}(z^{(i)})<-\gamma$, we think $g(z^{(i)})<0$. And $(-\gamma,\gamma)$ is the suspicious domain. While the value  $\hat{g}(z^{(i)})\in(-\gamma,\gamma)$, we calculate $g(z^{(i)})$ to verify. So the probability of $g(Z)<0$ can be represented as
Following the idea of the response surface method (RSM) \cite{au2001estimation,khuri2010response} and  design from \cite{li2011efficient}, the procedure of failure probability estimation can be specified as follows.
\begin{mdefinition}
	For a given real parameter $\gamma$, a limit state function $g$ and its surrogate model $\hat{g}$, the failure probability of $g(Z)<0$ can be represented as
	\begin{equation}\label{hy1}
	\begin{aligned}
	\mathrm{Prob}(g<0) \approx \mathrm{Prob} (\gamma;g,\hat{g}) = \mathrm{Prob}(\{\hat{g} < -\gamma\}) + \mathrm{Prob}(\{|\hat{g}| \leq \gamma\} \cap \{g < 0\}), 
	\end{aligned}
	\end{equation}
	where $(-\gamma,\gamma)$ is called the suspicious domain.
\end{mdefinition}
With MC method, we have
\begin{equation}\label{hy2}
\begin{aligned}
\widehat{\mathrm{Prob}}(g<0) &\approx \widehat{\mathrm{Prob}}(\gamma;g,\hat{g}) = \frac{1}{M}\sum_{i=1}^{M}(\chi_{\{\hat{g}<0\}}^{(h)}(z^{(i)}))\\
&=\frac{1}{M}\sum_{i=1}^{M}[\chi_{\{\hat{g}<-\gamma\}}^{(h)}(z^{(i)}) + \chi_{\{|\hat{g}| \leq \gamma\}}^{(h)}(z^{(i)}) \cdot \chi_{\{g<0\}}(z^{(i)})];
\end{aligned}
\end{equation}
where $\gamma$ in the equation is an arbitrary positive number, which denotes the range of re-verifying domain. It is clear that the computational cost will increase with $\gamma$ growing. So choosing a small enough $\gamma$ without sacrificing accuracy is a traditional difficulty. 

In this paper, we employ neural networks as the surrogate model $\hat{g}$, and we call the method as a direct combination of hybrid method and neural network as neural hybrid (NH) method.

\subsection{Multilayer neural network surrogate}
\begin{figure} 
	\centering
	\includegraphics[width=\linewidth]{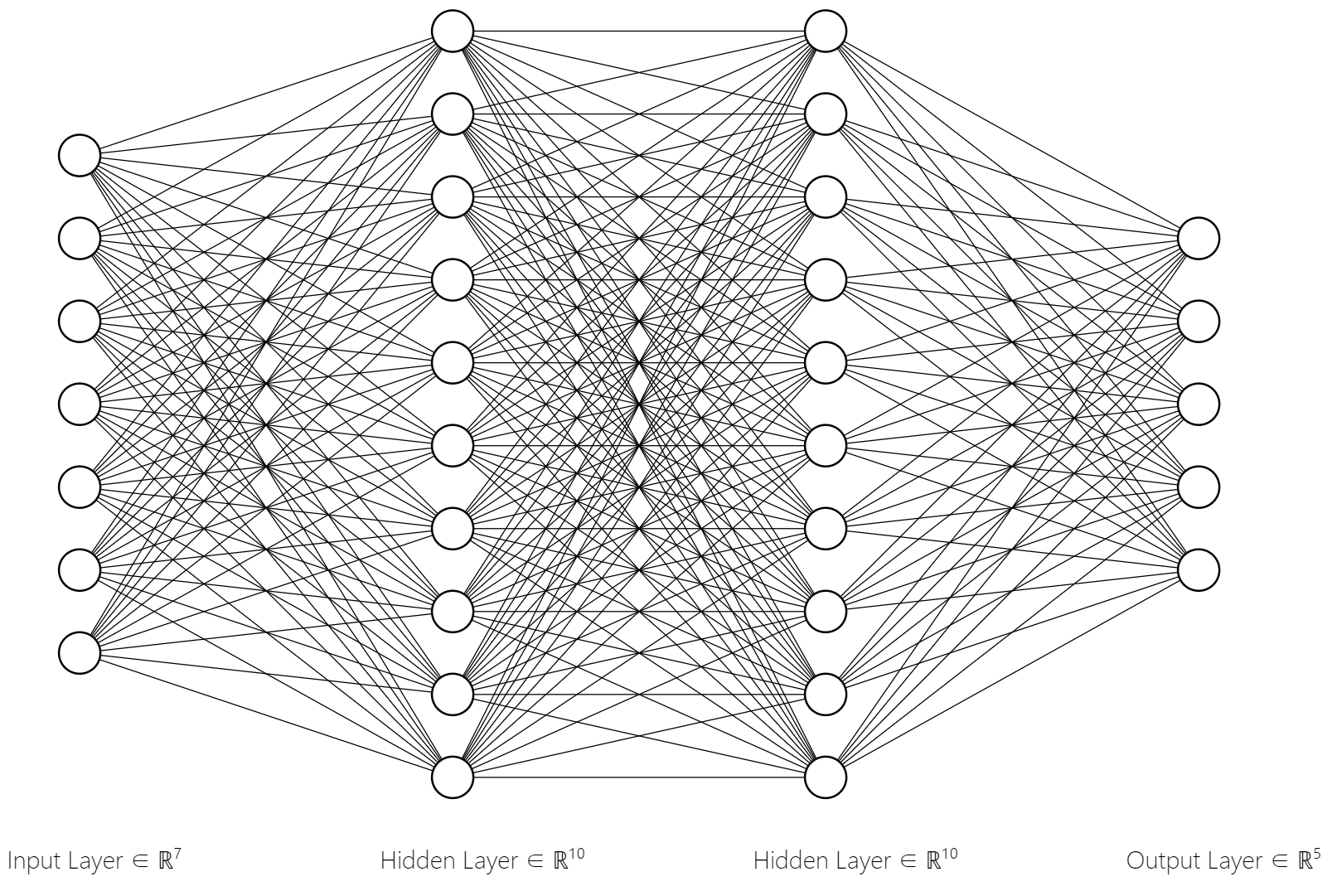}
	\caption{Schematic representation of a multilayer neural network which contains 7 input units, 5 output units and 2 hidden layers.}
	\label{ANN}
\end{figure}

 We employ full connected multilayer perceptron as the basic of our surrogate, which is a feed-forward type network with only adjoining layers connected. Network generally consists of input, hidden and output layers, see Figure \ref{ANN}. For illustration only, the network depicted consists of 2 layers with 10 neurons in each layer and $\sigma$ denotes an element-wise operator
 \begin{equation}
 \sigma(x)=\left(\phi\left(x_{1}\right), \phi\left(x_{2}\right), \ldots, \phi\left(x_{10}\right)\right),
 \end{equation}
 where $\phi$ is called activation function. We recall the network space with given data according to \cite{hornik1991approximation}
 \begin{equation}
 \mathfrak{N}_i^n(\sigma)=\left\{h(x) : \mathbb{R}^{k} \rightarrow \mathbb{R} | h(x)=\sum_{j=1}^{n} \beta_{j} \sigma\left(\alpha_{j}^{T} x-\theta_{j}\right)\right\}
 \end{equation}
 where $\sigma$ is any activation function, $\mb{x}\in\mathbb{R}^{k}$ is one set of observed data, $\beta\in\mathbb{R}^n$, $\alpha\in\mathbb{R}^{k\times n}$ and $\theta\in\mathbb{R}^{k\times n}$ denote coefficients of networks. The learning rule called back propagation algorithm \cite{rumelhart1986learning} is widely used in neural network. Multilayer perceptron is a mapping from $n-$dimension Euclidean space to $m-$dimension Euclidean space if there are n input units and m output units. {\ls We can also choose residual neural networks \cite{he2016deep} as surrogates for efficiency, 
 but we just present our idea of this work using the fully-connected neural network for simplicity.
 
 Via the comparison of the accuracy using neural networks with different layers in \cite{he2016deep}, the results of deeper neural networks are more accurate after sufficient training procedure. 
 Therefore,  we call deep networks fine-fidelity and shadow networks coarse fidelity in the following. In order to reduce the risk of overfitting, we employ cross validation and dropout \cite{srivastava2014dropout}.}

\section{Hierarchical neural hybrid method}
% Section 4

\begin{figure*}
	%\vskip -0.1in
	\centering{
		\includegraphics[width=0.8\textwidth,height=0.3\textheight]{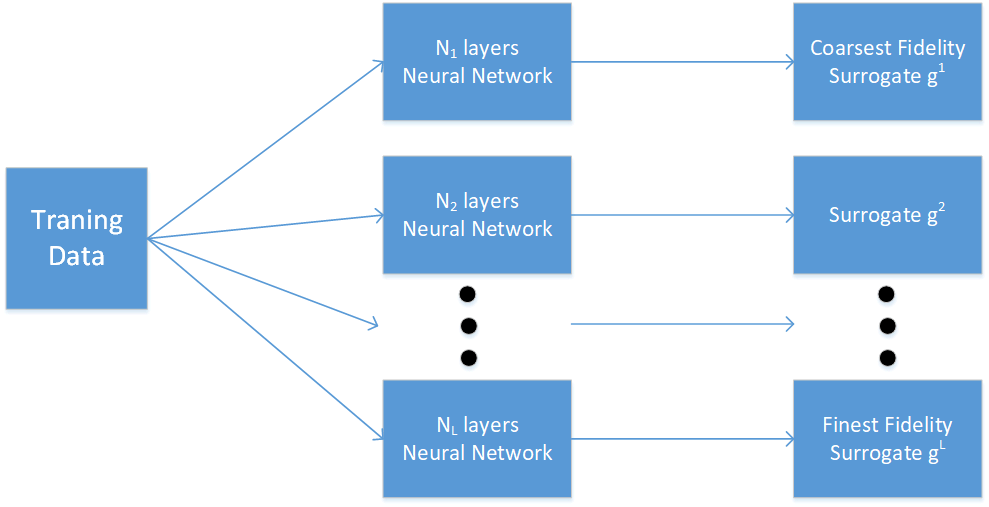}}
	\caption{Illustration of hierarchical training procedure.}
	\label{multi}
	\vskip -0.1in
\end{figure*}

In the former part, we introduced the hybrid method to speed up solving complex system by using an accurate enough surrogate. However, the computational cost of using a fine-fidelity model is also expensive, though it is much cheaper than the original system. Here we propose a hierarchical neural hybrid (HNH) method which constructs a hierarchy surrogate model to accelerate the standard NH method. The hierarchical surrogate models here are neural networks with different layers. Let  $g^{(1)},g^{(2)},\cdots,g^{(L)}$ denote $L$ surrogates, with $\ell=1$ being the coarsest and $\ell=L$ being the finest, each $g^{(\ell)} (\ell=1,\ldots,L)$ is a feedforward neural network with $P_\ell$ layers  and can approximate the limit state function $g$ well, illustrated in Figure \ref{multi}. It should be noticed that training data are the same for constructing hierarchical models off-line. Compared with NH method with single fine-fidelity, the hierarchical surrogate models can reduce the running time. 

\begin{figure*}
	\centering
%	\begin{tikzpicture}[scale=0.8, transform shape,node distance = 5.2cm, auto]
%	\node [block] (snn) {shallow neural network};
%	\node [block, right of=snn] (dnn) {deep neural network};
%	\node [blockplus, right of=dnn] (mltf) {multi-fidelity, cheap and accurate solution};
%	%\node [space, right of=dnn] (spacetwo) {};
%	\node [block, left of=snn] (true) {true model};
%	\node [decision, below of=snn] (lf) {biased solution};
%	\node [decision, below of=true] (exth) {groundtruth};
%	\node [cloud, above of=snn] (ours) {HNH (ours), cheap and optimal solution};
%	\node [decision, below of=dnn] (dnnf) {accurate solution};
%	
%	\path [line] (snn) -- node {coarse-fidelity, cheap}(lf);
%	\path [line] (dnn) -- node {fine-fidelity, expensive}(dnnf);
%	\path [line] (true) -- node {extremely expensive}(exth);
%	\draw [->,dashed] (snn) to [in = 135, out = 35] (mltf);
%	\draw [->, dashed] (snn) to (ours);
%	\draw [->, dashed] (dnn) to (ours);
%	\draw [->, dashed] (true) to (ours);
%	\draw [->, dashed] (dnn) to (mltf);
%	\end{tikzpicture}
	\includegraphics[width=0.8\textwidth]{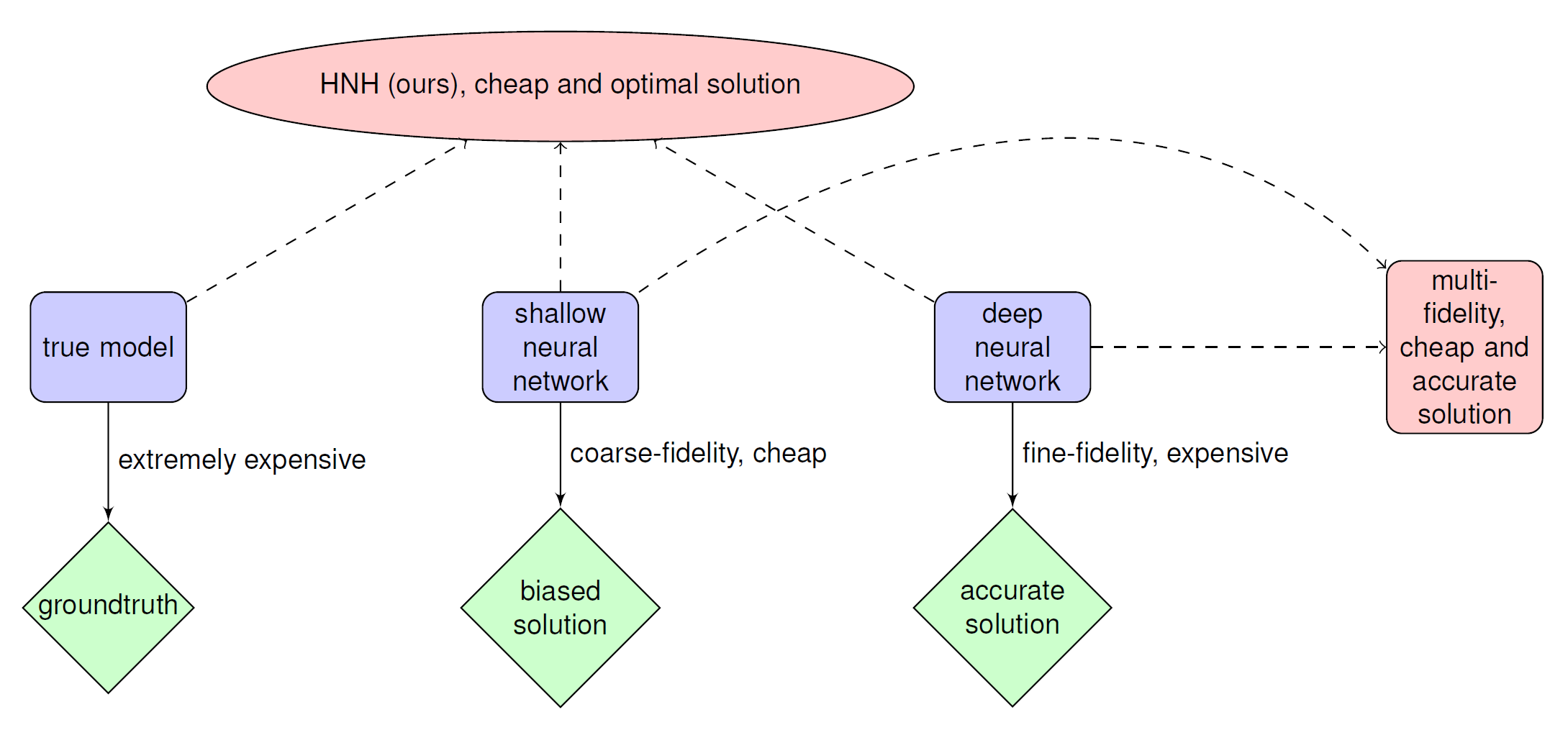}
	\caption{The idea of our HNH method.} 
	\label{form}
	\vskip -0.1in 
\end{figure*}

The idea of our method can be illustrated in Figure \ref{form}. 
{ Discrete solutions from true PDE models are referred as groundtruth, but the computational procedure can be extremely expensive. In former part we introduce that solutions from deep networks are more accurate but more expensive. If we use multifidelity to combine different networks, the cost can be reduced and the accuracy can be kept. Our hierarchical neural hybrid method uses the true model to improve the accuracy and the cost does not increase much, because our method only needs to solve 
the discrete PDEs a few times.} 
%Our HNH method is introduced in Section \ref{introHNH}-\ref{complex}. Section \ref{errorHNH} presents an analysis of the error and expectation, and we prove that our method is unbiased.

%\tikzstyle{decision} = [diamond, draw, fill=green!20, 
%text width=5.3em, text badly centered, node distance=3cm, inner sep=0pt]
%\tikzstyle{block} = [rectangle, draw, fill=blue!20, 
%text width=5em, text centered, rounded corners, minimum height=4em]
%%\tikzstyle{space} = [rectangle, draw, fill=white, 
%%text width=0, text centered, rounded corners, minimum height=0]
%\tikzstyle{line} = [draw, -latex']
%\tikzstyle{cloud} = [draw, ellipse,fill=red!20, node distance=3cm,
%minimum height=4em]
%\tikzstyle{blockplus} = [rectangle, draw, fill=red!20, 
%text width=5em, text centered, rounded corners, minimum height=4em]

\subsection{The HNH method}\label{introHNH}
\begin{figure}
	%\vskip -0.1in
	\centering{
	\includegraphics[width=0.8\textwidth]{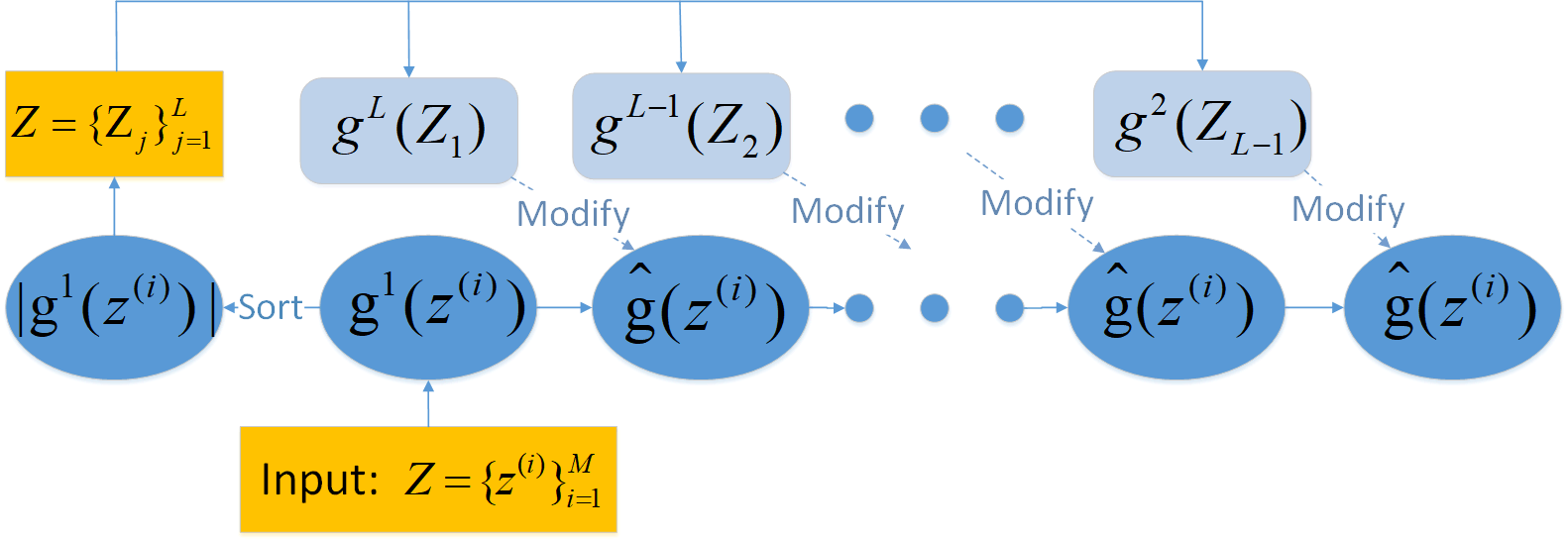}}
	\caption{Illustration of hierarchical neural hybrid method.} 
	\label{diag}
	\vskip -0.1in
\end{figure}

Hybrid method combines the robustness of MC and the feasibility of RSM. However, the computation cost is still expensive if an enough accurate surrogate model is repeatedly applied. Following the idea of multifidelity approaches {\ls\cite{ng2014multifidelity,aydin2019general,zhu2017multi,meng2019composite}}, we suppose that the cost can be reduced without losing accuracy by using a hierarchy surrogate model in the neural hybrid method. Our HNH method iterates through the level $\ell=0,\cdots,L$, where with $\ell>0$ increasing the accuracy increases except $g^{(0)}=g$ means the original system. First, the HNH method evaluates the state function at $Z$ by using $g^{(1)}$ instead of neural hybrid method using $g^{(L)}$, and then sorts $\{|g_1(z^{(i)})|\}_{i=1}^M$ in an ascending order. Second, the HNH method divides the data set into $L$ parts according to the value of $\{|g_1(z^{(i)})|\}_{i=1}^M$. For part $1,\cdots,L-1$, HNH employs $g^{(L)},\cdots,g^{(2)}$ to modify the prediction of failure with a threshold $\epsilon_{opt}$ to stop the modification. The modification procedure is presented in Algorithm \ref{HNH}. Third, our HNH method uses the modified results of prediction to run the iterative hybrid method introduced in \ref{IterAl} and Algorithm \ref{NN_Iter}. { A detailed illustration of HNH method can be found in Figure \ref{diag}.} 

\begin{algorithm}
	\caption{Modification procedure of hierarchical neural hybrid method}
	\label{HNH}
	\begin{algorithmic}[1]
		\State {\bfseries Input:} $g^{(1)},g^{(2)},\cdots,g^{(L)}$, $Z = \{z^{(i)}\}_{i=1}^{M}$, $\eta$;
		%\State {\bfseries Initialize:} $\epsilon=1$;
		\State Evaluate $g^{(1)}(z^{(1)}),g^{(1)}(z^{(2)}),\cdots,g^{(1)}(z^{(M)})$ and get label $\{\chi_{\{\hat{g}<0\}}(z^{(i)})\}_{i=1}^M = \{\chi_{\{g_1<0\}}(z^{(i)})\}_{i=1}^M$;
		\State Sort $\{z^{(i)}\}_{i=1}^M$ such that the value $\{|g^{(1)}(z^{(i)})|\}_{i=1}^M$ is in an ascending order, and the sorted order is denoted as $\{z^{(j)}\}_{j=1}^M$;
		\State Split dataset $Z=\{z^{(j)}\}_{j=1}^M$ into $L$ parts and each part has the same number of samples;
		\For {$i = 1$ {\bfseries to} $L-1$}
		\State $\epsilon=0$
		\For { $j = \frac{(i-1)M}{L}+1$ {\bfseries to} $\frac{iM}{L}$}
		\State $\Delta = [-\chi_{\{g^{(1)}<0\}}z^{(j)} + \chi_{\{g^{(-i+L-1)}<0\}}z^{(j)}]$
		\State $\chi_{\{\hat{g}<0\}}(z^{(j)}) = \chi_{\{g^{(1)}<0\}}(z^{(j)}) + \Delta$;
		\State $\epsilon = \epsilon + \frac{L}{M}[-\chi_{\{g^{(1)}<0\}}z^{(j)} + \chi_{\{g^{(-i+L-1)}<0\}}z^{(j)}]$;
		%\State $t = t + 1$;
		\EndFor
		%\State $\epsilon = \epsilon_{t}$
		\If {$\epsilon < \eta$}
		\State {\bfseries break}
		\EndIf
		\EndFor
		\State {\bfseries Return:} Label $\{\chi_{\{\hat{g}<0\}}(z^{(j)})\}_{j=1}^M$.
	\end{algorithmic}
\end{algorithm}

\subsection{Iteration algorithm of HNH method}\label{IterAl}
Here, we provide our integrated iteration algorithm for computing failure probability of high-dimensional problems through HNH method.

The hybrid method is applied which utilizes surrogate models to enhance the performance of Monte Carlo sampling. We present the probability formula defined in equation (\ref{hy1}) and (\ref{hy2}).

Hybrid method with $\gamma$ can enhance the performance of MC, but it is often hard to get the threshold $\gamma$. So we use iterative scheme to avoid choosing $\gamma$.

Let $\delta M$ be the number of samples generated from $g$ in each iteration, $\epsilon_{opt}$ be a given tolerance. The formal description of HNH method is presented in Algorithm \ref{NN_Iter}.
\begin{algorithm}
	\caption{Iteration algorithm of hierarchical neural hybrid method}
	\label{NN_Iter}
	\begin{algorithmic}[1]
		\State {\bfseries Input:} $g^{(1)},\cdots,g^{(L)}$, $Z = \{z^{(i)}\}_{i=1}^{M}$, $\eta$, $\delta M$, $\epsilon_{opt}$;
		\State {\bfseries Initialize:} $Z^{(0)}=\emptyset$, $k=0$, $\epsilon = 10 \epsilon_{opt}$;
		\State Obtain $\{\chi_{\{\hat{g}<0\}}(z^{(j)})\}_{j=1}^M$ using algorithm 1;
		\State $P_f^{(k)} = \frac{1}{M}\sum_{j=1}^{M}\chi_{\{\hat{g}<0\}}(z^{(j)})$;
		\While {$\epsilon > \epsilon_{opt}$}
		\State $k = k + 1$;
		\State $\delta Z^{(k)} = \{z^{(j)}\}_{j=(k-1)\delta M+1 }^{k\delta M}$;
		\State $Z^{(k)} =  Z^{(k-1)} \cup \delta Z^{(k)}$;
		\State $\delta P=\frac{1}{M}\sum\limits_{z^{(j)}\in\delta Z^{(k)}}[-\chi_{\{\hat{g}<0\}}(z^{(j)})+\chi_{\{g<0\}}(z^{(j)})]$;
		\State Update the failure probability: \par
		$P_f^{(k)} = P_f^{(k-1)} + \delta P$;
		\State $\epsilon = |P_f^{(k)} - P_f^{(k-1)}|$;
		\EndWhile
		\State {\bfseries Return:} $P_f^{(k)}$.
	\end{algorithmic}
\end{algorithm}

\subsection{The computational procedure and the computational complexity of HNH}\label{complex}
%Input models are $g^{(1)},g^{(2)},\cdots,g^{(L)}$, $Z$ is input data, $\eta$ is the threshold and $t$ means iteration times. 
We construct a hierarchy of neural network model denoted as $g^{(1)}$,$g^{(2)}$,$\cdots$,$g^{(L)}$, where the number of layers are $P_1,P_2,\cdots,P_L$ in a ascending order and there are $N$ neurons in each layer (see Algorithm \ref{HNH} ). Our input data set is $Z\in \mathbb{R}^{M\times r}, r\ll M$. For the classical single-fidelity neural network $g^{(L)}$, the computational complexity is $(M\times P_L\times N^2)$. In our HNH method, the computational complex is $(M\times P_1\times N^2 + \frac{M}{L}P_{L}N^2 + \cdots + \frac{M}{L}P_{L+1-\xi}N^2)$ if the modifications finished after $m$-th iteration, where $\xi = \mathrm{ceil}(\frac{mL}{M})$. The computational cost of our HNH method will be reduced significantly compared to NH method since $m\ll M$ and $\xi \ll L$ for an acceptable surrogate model.

\subsection{Analyses of HNH method}\label{errorHNH}
%\textbf{Assumption 1.} The models $g^{(\ell)}$ satisfy
%\begin{equation}
%|f^{(\ell+1)}(z)-f^{(\ell)}(z)|\leq \eta^{\ell};
%\end{equation}
%where $\eta\in(0,1)$, $\ell = 1,\cdots,L$.

\begin{massumption}\label{assum1}
Let $C>1$, $a>1$, $0<\rho<1$, $M$ is the size of input data $Z$, $\mathcal{F^{\ell}} = \{z\in Z | \{g^{(\ell)}(z)>\eta\}\cap\{g(z)<0\}\}$. The models $g^{(\ell)}$ satisfy
\begin{equation}
\int_{\Omega_{\eta}}\chi_{\mathcal{F^{\ell}}}dF_Z(z)\leq (\frac{1}{C}\frac{1}{1+\exp(ax)})^{\ell},
\end{equation}
%and
%\begin{equation}
%\eta_{max} = \Vert g-g^{(\ell)}\Vert_{\mathcal{L}_{\Omega}^p},
%\end{equation}
where $\mathcal{G}$ is sorted $\{|g^{(1)}(z^{(i)})|\}_{i=1}^M$, $\eta = \mathcal{G}_{\ell M/L}$, $\eta_{max} = \mathcal{G}_{(1-\rho)M}$, $x = \eta/\eta_{max}$.
\end{massumption}

This assumption takes the formula of inverse of the Sigmoid function. It is a reasonable assumption for it fitting the numerical results well.

\begin{massumption}\label{assum2}
For any $\epsilon>0$, given $C$ and $a$, $\eta_{max}$, $\eta_t$ is a threshold written as 
\begin{equation}
\eta_t = \ln(\frac{1}{C\cdot\epsilon})\frac{\eta_{max}}{a}.
\end{equation} 
\end{massumption}

\begin{mtheorem}
Failure probability is $P_f$, $P_f^h$ is calculated by neural hybrid method, $\hat{P}_f$ is the failure probability evaluated by HNH,  $g(Z)\in L_{\Omega}^p,p\geq1$ is the given system function, $\hat{g}(Z)$ is the HNH surrogate in $L^p-$ norm, $g^{(\ell)}(Z)$ is a hierarchy of surrogate models for $\ell = 1,\cdots,L$, for any $\epsilon>0$, there exists $\eta_t$, for any $\eta > \eta_t$,
\begin{equation}
|P_f - \hat{P}_f|\leq\epsilon.
\end{equation}
$Proof.$ 
\begin{equation}
|P_f - \hat{P}_f|\leq|P_f - P_f^h|+|P_f^h - \hat{P}_f|.	
\end{equation}
For $\ell = 1,\cdots,L$,
\begin{equation}
	\begin{aligned}\label{lammaP1}
		|P_f^h - \hat{P}_f| &= \sum_{\ell=1,\cdots,L-1}\int_{\Omega_{\eta_{\ell}}}\chi_{\{g^{L-\ell}(z)>\eta\}\cap\{g^{L}(z)<0\}}dF_z(z)\\
		&= \sum_{\ell=1,\cdots,L-1}\int_{\Omega_{\eta_{\ell}}}\chi_{\{\mathcal{U}-\{g^{L-\ell}(z)>\eta\}^{\complement}\cup\{g^{L}(z)<0\}^{\complement}\}}dF_z(z)
	\end{aligned}
\end{equation} 
where $\mathcal{U}$ means the universal set, $\eta_1 > \cdots > \eta_{L-1}$.

Choose $\ell$, for any $\epsilon>0$, there exists $\eta_t$, for any $\eta > \eta_t$, we can obtain
\begin{equation}\label{lammaP2}
	\begin{aligned}
		& \int_{\Omega_{\eta_{L-\ell}}}\chi_{\{g^{\ell}(z)>\eta\}\cap\{g^{L}(z)<0\}}dF_z(z)\\
		&= \int_{\Omega_{\eta_{L-\ell}}}\chi_{\{\mathcal{U}-\{g^{\ell}(z)>\eta\}^{\complement}\cup\{g^{L}(z)<0\}^{\complement}\}}dF_z(z)\\
		&= 1-(1-\frac{1}{C}\frac{1}{1+\exp(ax)})^{\ell} - (\frac{1}{C}\frac{1}{1+\exp(ax)})^{L}\\
		&\leq 1-(1-\frac{1}{C}\frac{1}{1+\exp(ax)})^{L} - (\frac{1}{C}\frac{1}{1+\exp(ax)})^{L}\\
		&\leq (\frac{1}{C}\frac{1}{1+\exp(ax)}) - (\frac{1}{C}\frac{1}{1+\exp(ax)})^{L} < \epsilon.
	\end{aligned}	
\end{equation}
Upon combining (\ref{lammaP1})-(\ref{lammaP2}), we obtain
\begin{equation}\label{lamma1}
	|P_f^h - \hat{P}_f| < \epsilon
\end{equation}

Then, with the definition of failure probability, we have
%\begin{equation}
%	\begin{aligned}
%		&P_f = \int_{\Omega_{\eta}}\chi_{\{g^{(L)}(Z)<-\eta\}\cap\{g(Z)<0\}}dF_Z(z) \\&+ \int_{\Omega_{\eta}}\chi_{\{|g^{(L)}(Z)|<\eta\}\cap\{g(Z)<0\}}dF_Z(z) \\&+ \int_{\Omega_{\eta}}\chi_{\{g^{(L)}(Z)>\eta\}\cap\{g(Z)<0\}}dF_Z(z).
%	\end{aligned}
%\end{equation}
%While $\eta > \eta_t$, for any $\epsilon$, we can obtain
\begin{equation}
	\begin{aligned}
		|P_f - P_f^h| 
		&= \int_{\Omega_{\eta}}\chi_{\{g^{(L)}(Z)>\eta\}\cap\{g(Z)<0\}}dF_Z(z) \leq (\frac{1}{C}\frac{1}{1+\exp(ax)})^{L} < \epsilon.
	\end{aligned}
\end{equation}
So that
\begin{equation}
	|P_f-\hat{P}_f|<\epsilon.
\end{equation}
\end{mtheorem}

%\subsection{HNH variance estimation}\label{ANOVA}
%Dataset $Z = \{z^{(i)}\}_{i=1}^M$, the mean and variance of samples are given by
%\begin{equation}
%    \mathbb{E}[g(Z)]=\frac{1}{M}\sum\limit_{i=1}^{M}g(z^{(i)}),
%\end{equation}
%and
%\begin{equation}
%    \mathbb{V}ar[g(Z)] = \frac{1}{M-1}\sum\limit_{i=1}^{M}(g(z^{(i)}-\mathbb{E}[g(Z)]))^2.
%\end{equation}
%Here, let $\hat{V}^{(\ell)}$ denote $\mathbb{V}ar[g^{(\ell)}(Z)]$, where $g^{(\ell)}$ is one of the hierarchy neural network surrogate models.

Before the estimation of expectation, we present the definition of multifidelity surrogate $\hat{g}(z)$ obtained by HNH.
\begin{mdefinition}(HNH Surrogate)
For input data $Z = \{z^{(i)}\}_{i=1}^M$, $M$ is the number of total samples, $g^{(\ell)},\ell=1,\cdots,L$ is a hierarchy of surrogate models, $m$ means the times of modifications and $\xi = ceil(\frac{mL}{M})$, then $\hat{g}$ can be presented as
\begin{equation}\label{def}
\begin{aligned}
    \hat{g} &= \frac{1}{L}\sum\limits_{i=1}^{\xi-1}g^{(L+1-i)} + \frac{mL-(\xi-1)M}{ML}g^{(L-\xi+1)} + \frac{M-m}{M}g^{(1)},
\end{aligned}
\end{equation}
\end{mdefinition}
when $\xi=1$, the first term at the right side equals zero.

\begin{massumption}
$g(Z)$ is given system function, $g^{(\ell)}(Z)$ is a hierarchy of surrogate models for $\ell = 1,\cdots,L$, $Z$ is the dataset. While the models are trained with labeled data, the neural networks can be regarded as unbiased approximators. The expectations are given as
\begin{equation}
\mathbb{E}[g] = \mathbb{E}[g^{(1)}] = \cdots = \mathbb{E}[g^{{(L)}}].
\end{equation}
\end{massumption}

\begin{mtheorem}
The HNH surrogate $\hat{g}(Z)$ is an unbiased estimator of $g(Z)$. 
%\begin{equation}
%\mathbb{E}[\hat{g}] = \mathbb{E}[g].
%\end{equation}

\noindent$Proof.$ Suppose HNH modifications finish after $m$-th iteration, $\xi = ceil(\frac{mL}{M})$,
%\begin{equation}
%\begin{aligned}
%\hat{g} &= \frac{1}{M}[\frac{M}{L}[g^{(1)}+(-g^{(1)}+g^{(L)})] \\&+\cdots\\&+\frac{M}{L}g^{(1)}+\frac{mL-M\xi+M}{L}(-g^{(1)}+g^{(L-\xi+1)})\\&+\frac{(L-\xi)M}{L}g^{(1)}].
%\end{aligned}
%\end{equation}

\begin{equation}
\begin{aligned}
\mathbb{E}[\hat{g}] &= \frac{1}{M}[\frac{M}{L}\mathbb{E}[g^{(1)}+(-g^{(1)}+g^{(L)})] +\cdots +\frac{M}{L}\mathbb{E}[g^{(1)}]\\&+\frac{mL-M\xi+M}{L}\mathbb{E}[(-g^{(1)}+g^{(L-\xi+1)})] +\frac{(L-\xi)M}{L}\mathbb{E}[g^{(1)}]] \\&=\frac{1}{M}[\frac{M}{L}L\cdot\mathbb{E}[g]] = \mathbb{E}[g].
\end{aligned}
\end{equation}
\end{mtheorem}
%Therefore our multifidelity model $\hat{g}$ (\ref{def}) is an unbiased surrogate of $g$. This shows that $\hat{P_f}$ is an unbiased estimator of failure probability $P_f$.

%\begin{mproposition}
%Follow Definition 1, $M$ is the number of samples, let $g_i,i=1,\cdots,M$ present $M$ components of $\hat{g}$, 
%\begin{equation}
%    \mathbb{V}ar[\hat{g}] = \sum\limits_{i=1}^M\sum\limits_{j=1}^M\Cov(g_i,g_j).
%\end{equation}
%$Proof.$
%\begin{equation}
%    \begin{aligned}
%    &\mathbb{V}ar[\hat{g}] = \mathbb{E}[(\hat{g}-\mathbb{E}[\hat{g}])^2]
%    \\&=\mathbb{E}[(\sum\limits_{i=1}^Mg_i-\mathbb{E}(\sum\limits_{i=1}^Mg_i))(\sum\limits_{j=1}^Mg_j-\mathbb{E}(\sum\limits_{j=1}^Mg_j))]
%    \\&=\sum\limits_{i=1}^M\sum\limits_{j=1}^M\mathbb{E}[(g_i-\mathbb{E}[g])(g_j-\mathbb{E}[g])]=\sum\limits_{i=1}^M\sum\limits_{j=1}^M\Cov(g_i,g_j).
%    \end{aligned}
%\end{equation}
%\end{mproposition}

\section{Numerical experiments}
% Section 6: Numerical study
In this section we provide three numerical examples to test the performance of our HNH method. For the multivariate benchmark, we use MC to obtain the reference solution. In the test of diffusion and Helmholtz equations, we employ finite element method to calculate the solution. For these three studies, we trained hierarchical neural networks with 6, 15 and 30 layers, and each layer has 500 neurons. {\ls The selection of the number of layers and neurons is empirical.} All timings conduct on an Intel Core i5-7500, 16GB RAM, Nvidia GTX 1080Ti processor with MATLAB 2018a and Tensorflow under Python 3.6.5. The number of samples we generated in the following tests is related to the limit of RAM. In the time comparison, we set the runtime of HNH method as a unit. 

In numerical tests, the numbers of training samples are $10^3$, and the total cost includes these parts. While traditional approaches are expensive or infeasible in high-dimensional problems, we compare with the neural hybrid method proposed in this work.
{\ls It should be noted that after modifications using fine-fidelity networks, the only difference between HNH and NH is runtime, which means that the error plots of two methods in the following 
are the same with respect to
the number of discrete PDE solves.}

% We define the speedup rate as 
% \begin{equation}
% T_a = \frac{|t_o-t_n|}{t_n},
% \end{equation}
% where $T_a$ means time accelerated, $t_o$ means runtime using neural hybrid method, $t_n$ means runtime using HNH method.

\subsection{Multivariate benchmark}
Here we consider a high-dimensional multivariate benchmark problem in the field of structural safety from \cite{engelund1993benchmark}.
\begin{equation}
\begin{aligned}
g({ X}) &= \beta n^{\frac{1}{2}}-\sum_{i=1}^{n}X_i,
\end{aligned}
\end{equation} 
where $\beta =3.5$, $n=50$ and ${ X_i} \sim N(0,1)$.

Then we construct our surrogate model using neural network as
\begin{equation}\label{NN_surr}
\hat{g} = \phi_n(w_{n}\phi_{n-1}(\cdots\phi_1(w_1x+b_1)+\cdots)+b_{n}),
\end{equation}
where $\{\phi_j\}_{j=1}^n$ are projections, $\{w_j\}_{j=1}^n$ are weights, $\{b_j\}_{j=1}^n$ are biases and $x$ is input variable generated from the given distribution.

We consider the failure probability $P_f=\mathrm{Prob}(g(X)<0)$. With $5\times10^6$ MC sampling, the reference estimation is $P_f^{mc} = 2.218\times10^{-4}$.

% \begin{figure}[H]
% 	\begin{minipage}[t]{0.5\textwidth}
% 		\centering
% 		\includegraphics[width=\textwidth]{image/absolute_error.eps}
% 	\end{minipage}
% 	\begin{minipage}[t]{0.5\textwidth}
% 		\centering
% 		\includegraphics[width=\textwidth]{image/relatively_error.eps}
% 	\end{minipage}
%     \caption{Example 1, absolute error(left) and relatively error(right) of failure probability estimates for different iteration times .}
% \end{figure}

\begin{figure}[H]
\centering
\subfigure[Absolute error]{
		\centering
		\includegraphics[width=0.4\textwidth]{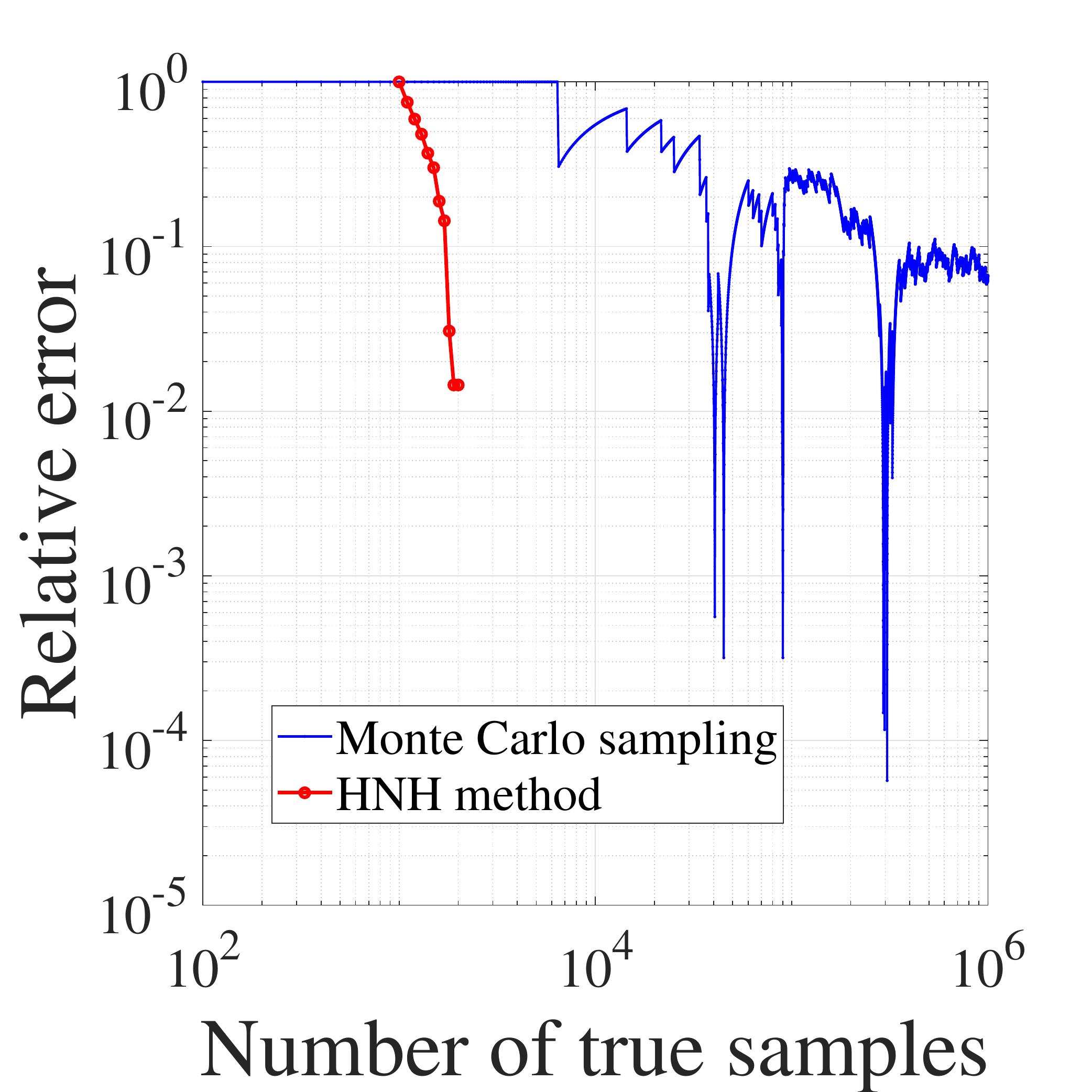}
		\label{2.1}}
%\hspace{0.001in}
\subfigure[Time comparison]{
		\centering
		\includegraphics[width=0.4\textwidth]{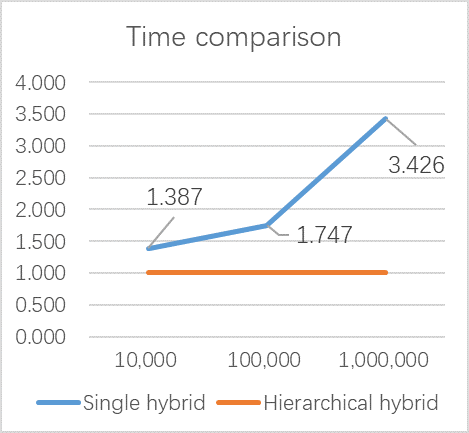}
		\label{2.2}}
	\caption{Performances on accuracy and time of HNH in multivariate benchmark.}
	\label{Fig_2}
\end{figure}

We generated $5\times10^6$ samples for evaluating reference estimation, so we set $10^6$ samples for the error estimation. In Figure \ref{2.1}, blue line is the error of Monte Carlo estimation, and red line is the error of probability evaluated by HNH method. The error decreases quickly and iteration stops automatically since the $\epsilon$ between the last two iterations is lower than the threshold $\epsilon_{opt}$ \textbf{(last two red circles are too close to identify in Figure \ref{2.1}). The phenomenon perfectly fits our Assumption \ref{assum1} and \ref{assum2}, for the error decreasing to zero when $g(z)\geq\eta_t$.} It means that there will be no more modification when $g(z)\geq\eta_t$, and the failure probability estimation converges. So an acceptable probability ($P_f=2.25\times10^{-4}$) can be obtained with $2\times10^3$ sample generated from origin system, where $10^3$ samples are for training models. Compared with reference value which obtained by using $5\times10^6$ MC sampling, the absolute error is about $10^{-6}$ and relative error is $1.443\%$ with less than $0.4\%$ training samples generated from origin system $g$. In Figure \ref{2.2}, we show the timing comparison between our HNH method and our NH method. We set the running time of NH method as the unit time in three orders of magnitude. With number of samples growing, our HNH method is more efficient.

\subsection{Diffusion equation}
In this numerical test, we consider a diffusion problem. The governing equations of the diffusion problem are
\begin{equation}\label{FEM}
\begin{aligned} - \nabla \cdot ( a ( x , \xi ) \nabla u ( x , \xi ) ) & = 1 \quad \text { in } \quad D \times \Gamma \\ u ( x , \xi ) & = 0 \quad \text { on } \quad \partial D _ { D } \times \Gamma \\  \frac { \partial u ( x , \xi ) } { \partial n } &= 0 \quad \text { on } \quad \partial D _ { N } \times \Gamma.\end{aligned} 
\end{equation}
Where the equation setting can be found in \cite{liao2016reduced}. We can employ finite element method (FEM) and the weak form of (\ref{FEM}) is to find an appropriate $u ( x , \xi ) \in H _ { 0 } ^ { 1 } ( D )$ s.t. $\forall v \in H _ { 0 } ^ { 1 } ( D ) , ( a \nabla u , \nabla v ) = ( 1 , v )$. And the finite element solver is from the former work \cite{liao2013implicit,elman2014ifiss}.

In our numerical study, the spatial domain is $D(0,1)\times(0,1)$. Dirichlet boundary conditions are applied on the left $(x=0)$ and right $(x=1)$ boundaries. Neumann conditions are applied on the top and bottom conditions. The problem is discretized on a uniform $65\times65$ grid, and $N_h = 4225$ is the spatial degrees of freedom.  

The coefficient $a ( x , \xi )$ of the diffusion problem is regarded as a random field, where $a_0(x)$ is mean function, $\sigma$ is standard deviation and covariance function,
\begin{equation}\label{cov}
\operatorname { Cov } ( x , y ) = \sigma ^ { 2 } \exp \left( - \frac { \left| x _ { 1 } - y _ { 1 } \right| } { L } - \frac { \left| x _ { 2 } - y _ { 2 } \right| } { L } \right),
\end{equation}
where $L$ is the correlation length. We employ Karhunen-Loéve (KL) expansion \cite{ghanem1991stochastic,babuvska2007stochastic} to approximate the random field
\begin{equation}
a ( x , \xi ) \approx a _ { 0 } ( x ) + \sum _ { k = 1 } ^ { d } \sqrt { \lambda _ { k } } a _ { k } ( x ) \xi _ { k },
\end{equation}
where $a _ { k } ( x )$ and $\lambda _ { k }$ are the eigenfunctions and eigenvalues of (\ref{cov}), $\{\xi_k\}_{k=1}^d$ are random variables. We set $a_0(x) = 1$, $\sigma = 0.42$, and the number $d$ we choose is large enough to capture 95\% of total variance of the exponential covariance function \cite{powell2009block}. In this paper, we set $d=48$ for correlation length $L=0.8$.

\begin{figure}
	\centering{
	\includegraphics[width=0.8\textwidth]{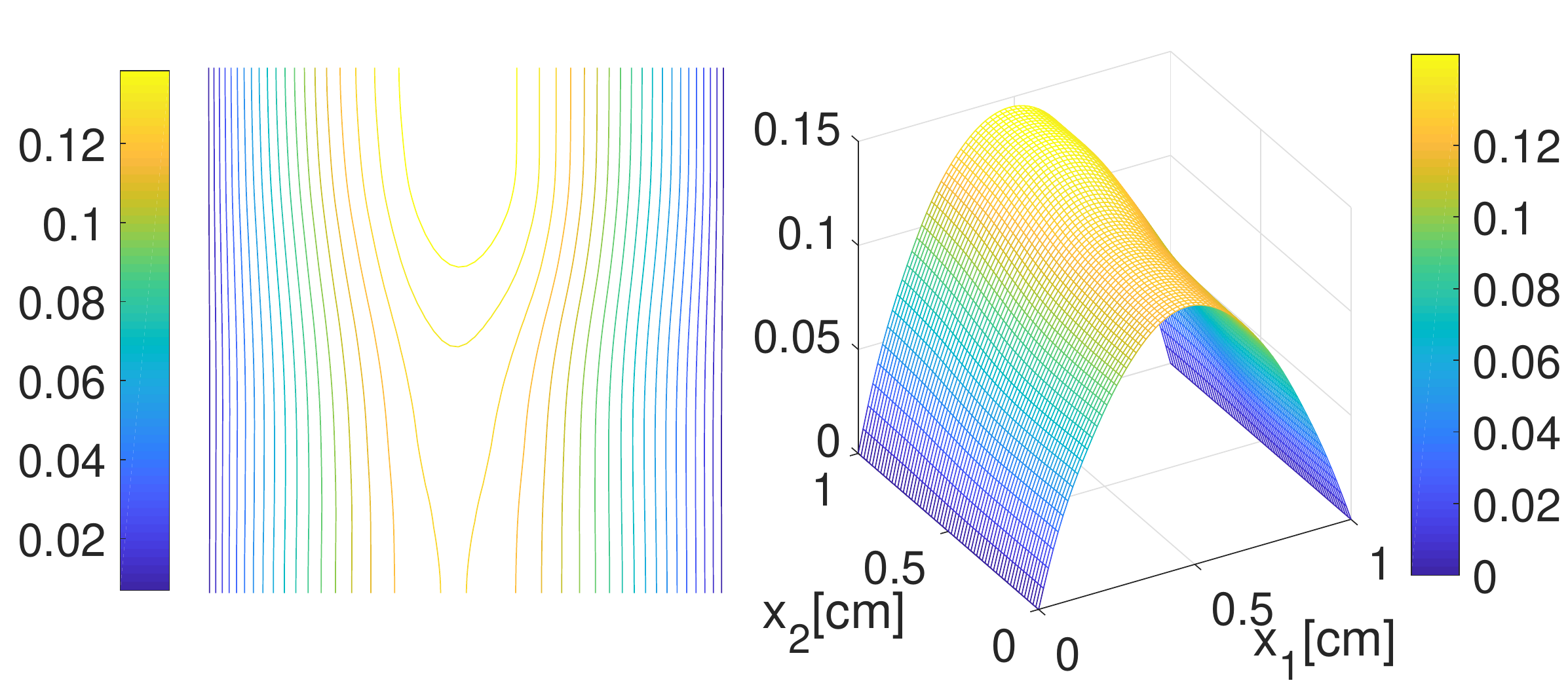}}
	\caption{The plots show the values of the diffusion equation, where left one is the contour.} 
	\label{Fig_3}
	\vskip -0.2in
\end{figure}

The numerical results solved by FEM using IFISS \cite{elman2014ifiss} are shown in Figure \ref{Fig_3}, and we choose $X = [0.5;0.5]$ as the point sensor placed. For $d=48$, with $10^6$ MC sampling, the reference failure probability for $\mathrm{Prob}(u(x,\xi)>0.19)$ is $P_f=1.2\times10^{-3}$.

\begin{figure}
\centering
\subfigure[][Absolute error]{
		\centering
		\includegraphics[width=0.4\textwidth]{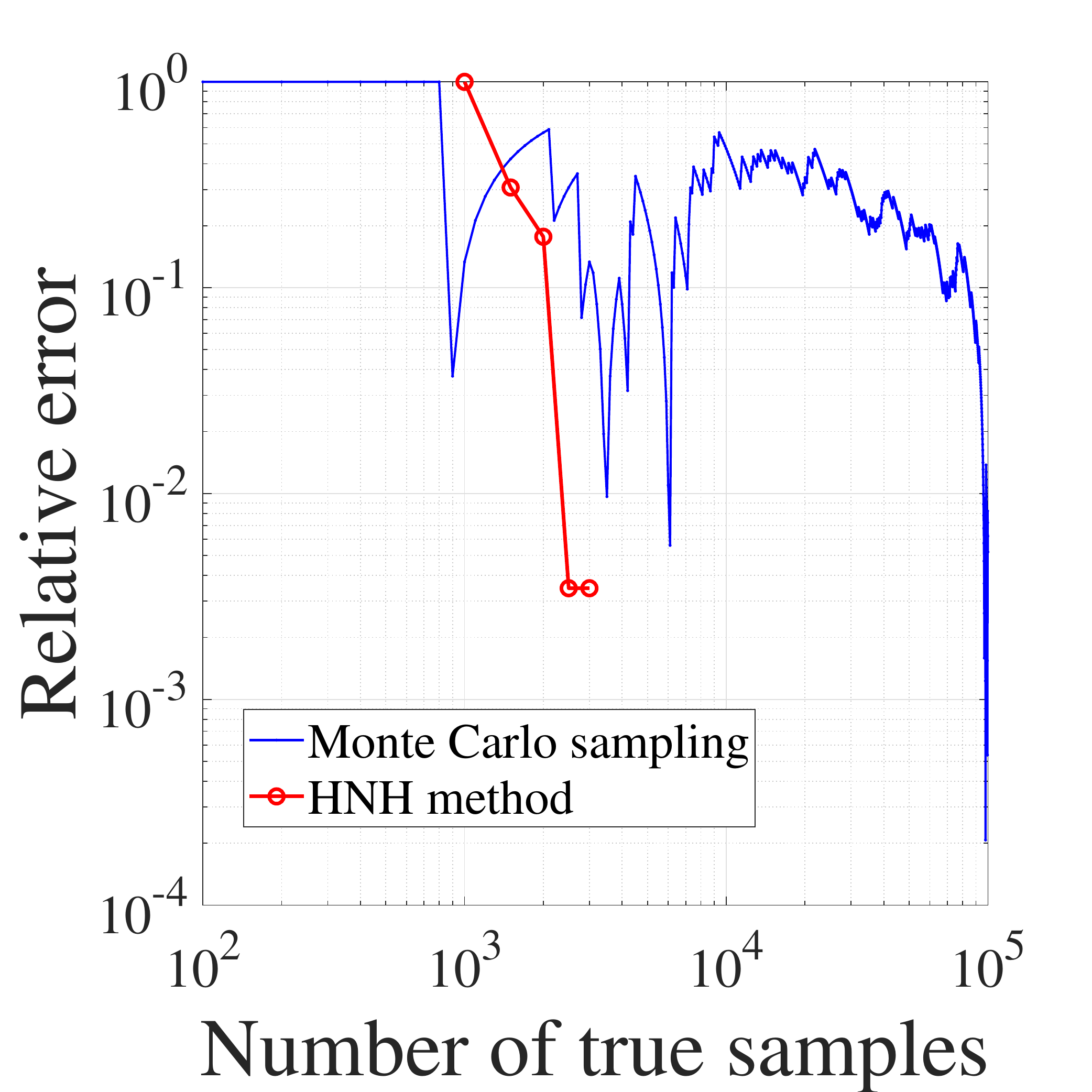}
		\label{4.1}}
\subfigure[][Time comparison]{
		\centering
		\includegraphics[width=0.4\textwidth]{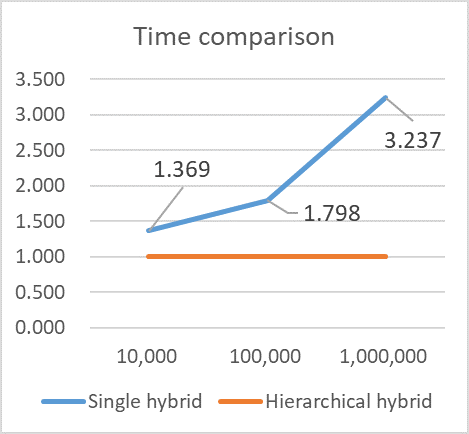}
		\label{4.2}}
	\caption{Performances on accuracy and time of HNH in diffusion equation for $d=48$.}
	\label{Fig_4}

\end{figure}

 In Figure \ref{4.1}, our HNH method uses $2\times10^3$ samples from origin model, and it captures the the failure probability as $1.15\times10^{-3}$, relative error is $4.17\%$ (the reference solution computed by MC using $10^6$ samples ). The result is better than MC with around $10^5$ samples. In Figure \ref{4.2}, our HNH method is more than three times faster than NH method when the magnitude is $10^6$. Surrogate method is widely used in huge magnitude problem. NH method is an efficient surrogate method, and our HNH method is more outstanding.

\subsection{Helmholtz equation}
We now consider the Helmholtz equation
\begin{equation}
- \Delta u - k ^ { 2 } u = 0,
\end{equation}
where $k$ is coefficient, and we set the homogeneous term.

\begin{figure}
\centering
	\begin{minipage}[t]{0.4\textwidth}
		\centering
		\includegraphics[width=\textwidth]{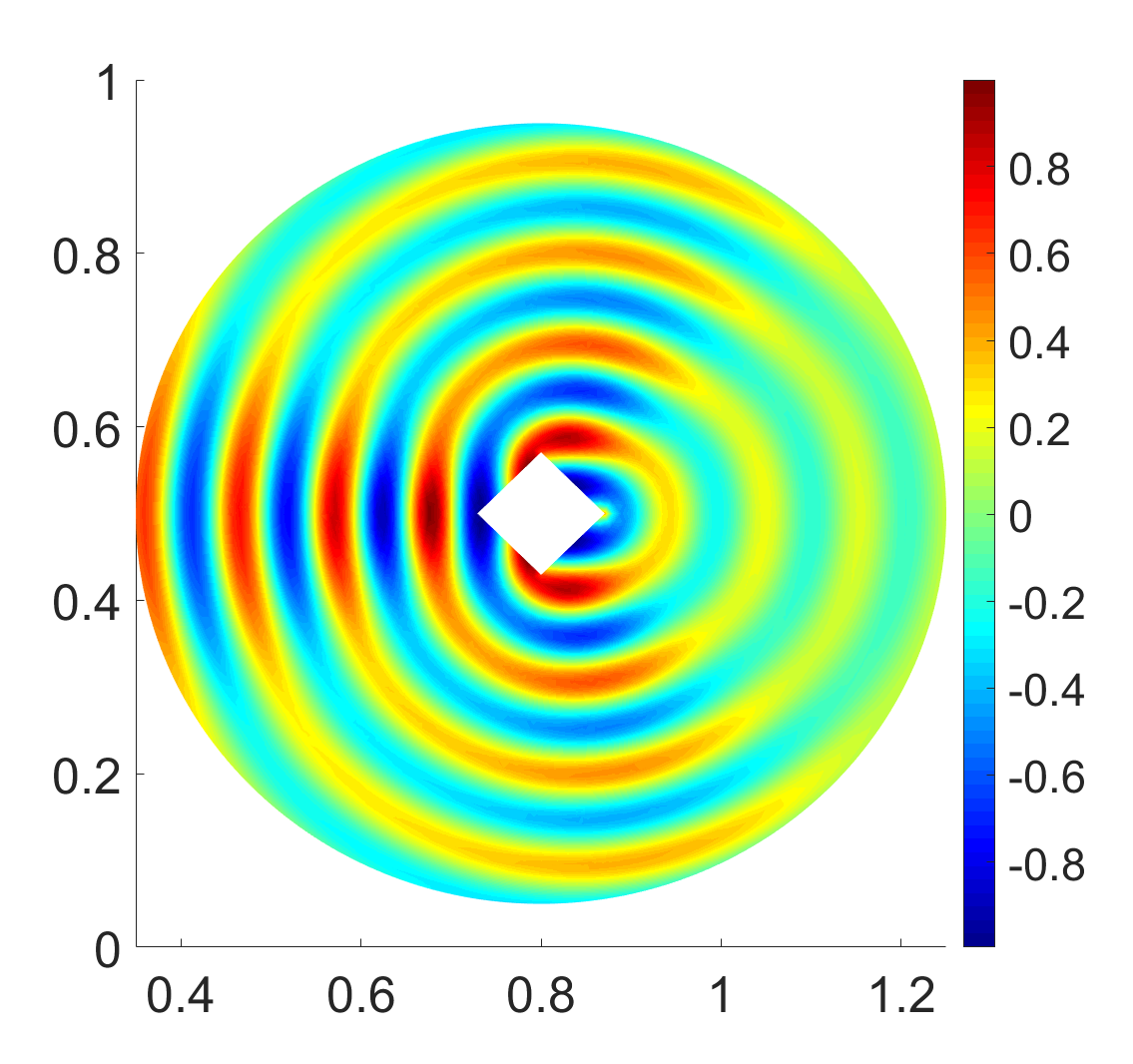}
		\label{5.1}
	\end{minipage}
	\begin{minipage}[t]{0.4\textwidth}
		\centering
		\includegraphics[width=\textwidth]{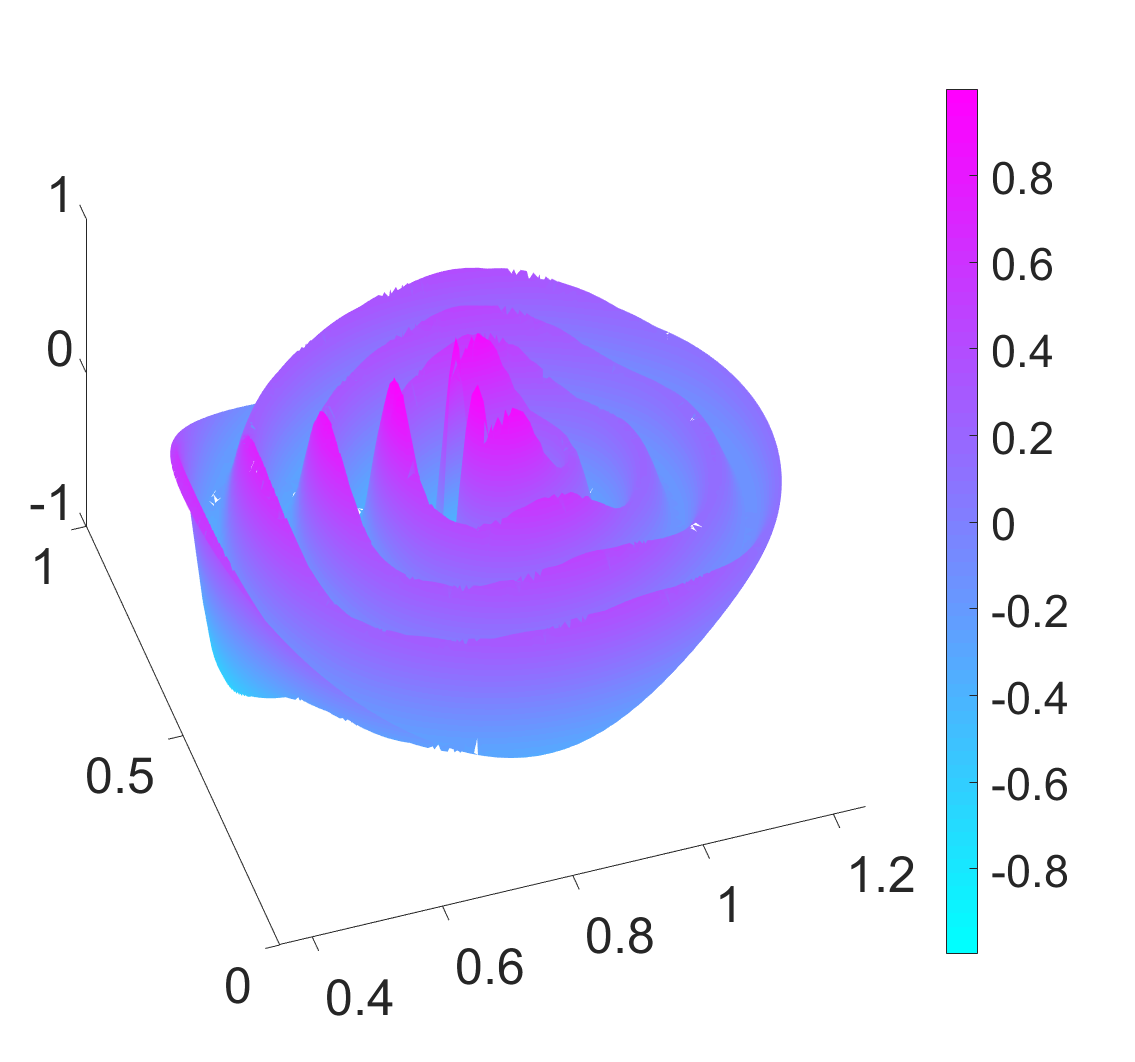}
		\label{5.2}
	\end{minipage}
	\caption{A snapshot of the value of the Helmholtz equation.}
	\label{Fig_5}

\end{figure}

We employ MATLAB PDE solver to obtain accurate solutions of Helmholtz equation shown in Figure \ref{Fig_5}. We choose $X = [0.7264;0.4912]$ as the point sensor placed. The reference failure probability for $\mathrm{Prob}(u(x,k)>1.09)$ is $P_f=2.08\times10^{-3}$ computed by MC using $10^5$ samples.

\begin{figure}
\centering
\subfigure[][Absolute error]{
		\centering
		\includegraphics[width=0.45\textwidth]{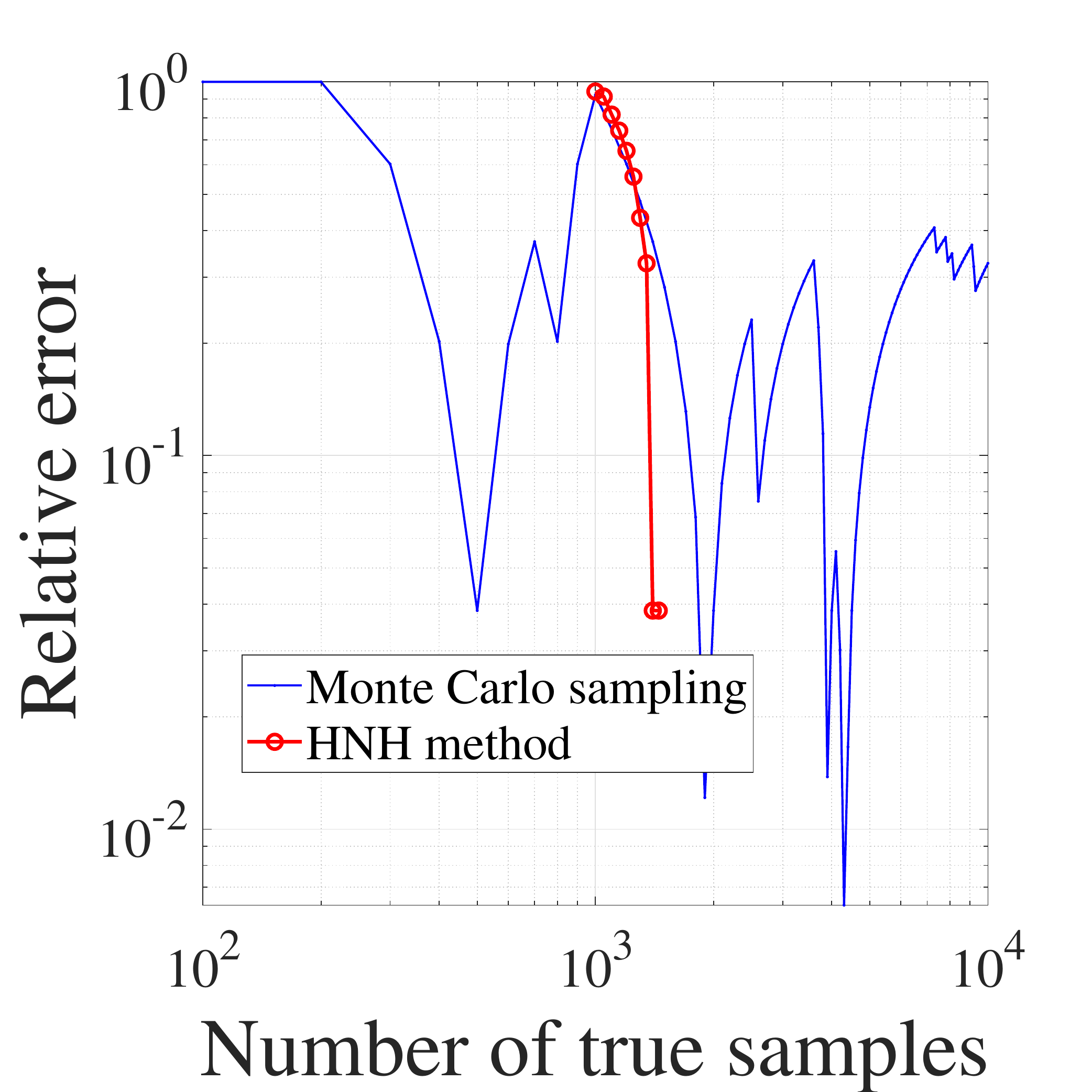}
		\label{6.1}}
\subfigure[][Time comparison]{
		\centering
		\includegraphics[width=0.45\textwidth]{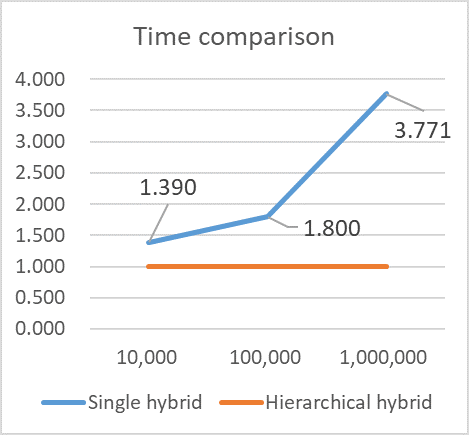}
		\label{6.2}}
	\caption{Performances on accuracy and time of HNH in Helmholtz equation.}
	\label{Fig_6}

\end{figure}

Figure \ref{6.1} shows the relative error compared with reference failure probability. Our method captures the failure probability ($P_f=2.16\times10^{-3}$) with only $1500$ samples generated from Helmholtz equation where $10^3$ for training and $500$ for modification. The absolute error is about $8\times10^{-5}$ and relative error is $3.85\%$. The result is more accurate than MC with $10^4$ samples. The estimation of standard MC does not converge within $10^4$ samples. Figure \ref{6.2} shows the speedup of our HNH approach compared with our former neural hybrid method. Considering the performance in both accuracy and efficiency, our HNH approach performs well.

% \begin{table}[htbp]
%  \caption{\label{tab:test} Time comparison}
%  \begin{tabular}{ccccc}
%   \toprule
%   Model & Size & Bench & Diff & Helm \\
%   \midrule
%  HNH & $10^4$ & 1 & 1 & 1 \\
%  NH & & 1.22 & 1.26 & 1.17 \\
%  HNH & $10^5$ & 1 & 1 & 1 \\
%  NH & & 1.18 & 1.27 & 1.31 \\
%  HNH & $10^6$ & 1 & 1 & 1 \\
%  NH & & 1.36 & 1.33 & 1.54 \\
%   \bottomrule
%  \end{tabular}
% \end{table}

% \begin{figure}[H]
% \vskip -0.1in
% \begin{minipage}[t]{0.45\textwidth}
% 		\centering
% 		\includegraphics[width=\textwidth]{image/diffcompar.eps}
% 		\label{7.1}
% 	\end{minipage}
% 	\caption{Time comparison for three numerical studies.}
% 	\label{Fig_7}
% \vskip -0.2in
% \end{figure}

\section{Conclusion}

Conducting adaptivity is a main concept for efficiently estimating failure probabilities of complex PDE models with high-dimensional inputs.
In this work, our HNH procedure adaptively fits the hierarchical structures for this problem. 
To finally show the efficiency  of HNH, we compare it with a direct combination of  neural network and hybrid method (which is referred to as NH).
Table \ref{table} shows the running times of HNH and NH to achieve the same given accuracy for the three test problems.
It can be seen that, as the sample size increases, the efficiency of HNH becomes more clear, e.g., for the Helmholtz problem with  $M=10^6$,
the running time of HNH is around a quarter of the time of NH.
Moreover, our method employs neural network as a surrogate to overcome the limitations of standard polynomial chaos for high-dimensional problems.
From the numerical examples, to achieve the same accuracy, it is clear that HNH only solves the PDEs several thousand times, while traditional MC 
needs to solve PDEs more than $10^5$ times.
Due to the universal nature of the neural network, our HNH method can be extend to general surrogate modelling problems.

\begin{table}[H]
	\centering
	\begin{threeparttable}
		\caption{Runtime for HNH method and its corresponding neural hybrid (NH) method with different number of samples for three examples in the former part.}
		\label{table}
		\begin{tabular}{cccc}
			\toprule 
			Type&$10^4$ / (ms)&$10^5$ / (ms)&$10^6$ / (ms)\\
			\midrule
			NH in benchmark&374.1&556.5&2659.5\\
			HNH in benchmark&269.6&318.6&776.4\\
			NH in diffusion&365.6&563.8&2630.7\\
			HNH in diffusion&267.1&313.5&812.7\\
			NH in Helmholtz&359.4&536.5&2449.3\\
			HNH in Helmholtz&258.5&298.1&649.5\\
			\bottomrule
		\end{tabular}
	\end{threeparttable}
\end{table}
\section*{Acknowledgments}
This work is supported by the National Natural Science Foundation of China (No. 11601329 and No. 11771289) and the science challenge project (No. TZ2018001). 
\newpage
\addcontentsline{toc}{chapter}{Bibliography}
\bibliography{KeLi_hnh}

\end{document}